\documentclass{article}
\usepackage{epsfig,epic,eepic,latexsym, amssymb}
\input xy

\xyoption{all}
\xyoption{all}
\makeindex

\def \pf{PROOF:}
\def \QED{\hfill\hbox{\hskip 4pt
                \vrule width 5pt height 6pt depth 1.5pt}}
\def \epf{\QED\\}

\newcommand{\Cl}{ \{   \hskip -3pt \mid  }
\newcommand{\Cr}{ \mid \hskip -3pt \}  }
\newcommand{\Sl}{ [ \hskip -1.5pt [  }
\newcommand{\Sr}{ ] \hskip -1.5pt ]  }

\newcommand{\V}{  {\cal V} }
\newcommand{\A}{  {\cal A} }
\newcommand{\B}{  {\cal B} }
\newcommand{\C}{  {\cal C} }
\newcommand{\D}{  {\cal D} }
\newcommand{\K}{  {\cal K} }
\newcommand{\LL}{  {\cal L} }

\newcommand{\Q}{ {\cal Q} }

\newcommand{\bicatB}{ {\bf B}}

\newcommand{\Vp}{ {\cal V}' }

\newcommand{\VCat}{ {\cal V}\mbox{-}\mathit{\bf Cat}}
\newcommand{\VCAT}{ {\cal V}\mbox{-}\mathit{\bf CAT}}
\newcommand{\VMod}{ {\cal V}\mbox{-}\mathit{\bf Mod}}
\newcommand{\VMOD}{ {\cal V}\mbox{-}\mathit{\bf MOD}}

\newcommand{\Cat}{ \mathit{\bf Cat}}
\newcommand{\CAT}{ \mathit{\bf CAT}}

\newcommand{\VpCat}{ {\cal V}'\mbox{-}\mathit{\bf Cat}}
\newcommand{\VpCAT}{ {\cal V}'\mbox{-}\mathit{\bf CAT}}
\newcommand{\SCat}{ Set\mbox{-}\mathit{\bf Cat} }
\newcommand{\SCAT}{ Set\mbox{-}\mathit{\bf CAT} }

\newcommand{\I}{\cal I}

\newcommand{\Co}{\mathit{\bf Cocts}}
\newcommand{\Cont}{\mathit{\bf Conts}}

\newcommand{\DCont}{ \Delta\mbox{-}\mathit{\bf Conts}}
\newcommand{\VCont}{ {\cal V}\mbox{-}\mathit{\bf Conts}}
\newcommand{\PhiCo}{\Phi\mbox{-}\mathit{\bf Cocts}}
\newcommand{\PhiCont}{\Phi\mbox{-}\mathit{\bf Conts}}

\newcommand{\PsiCo}{\Psi\mbox{-}\mathit{\bf Cocts}}
\newcommand{\PsiCont}{\Psi\mbox{-}\mathit{\bf Conts}}

\newcommand{\ZCont}{ 0\mbox{-}\mathit{\bf Conts}}

\newcommand{\p}{{\cal P}}
\newcommand{\PCo}{{\cal P}\mbox{-}\mathit{\bf Cocts}}
\newcommand{\PCont}{{\cal P}\mbox{-}\mathit{\bf Conts}}

\newcommand{\Rl}{rlift}
\newcommand{\Ll}{llift}

\newtheorem{theorem}{Theorem}[section]

\newtheorem{proposition}[theorem]{Proposition}
\newtheorem{lemma}[theorem]{Lemma}

\newtheorem{remarks}[theorem]{Remarks}
\newtheorem{remark}[theorem]{Remark}

\newtheorem{fact}{}[section] 

\title{Notes on enriched categories with colimits of some class}
\author{G.M. Kelly and V.Schmitt}

\begin{document}
\maketitle                      

\begin{abstract}
This work is on-going. Results presented here are issued mainly
from unpublished notes of the first author and contain those in
the notes http://arXiv.org/abs/math.CT/0309209 and
http://arXiv.org/abs/math.CT/0403164 of the second author.
The latter is also the only responsible for the typos, spelling
mistakes and other errors that could occur in the present paper!
We intend to add further sections 
to this work in the near future.\\

Given a class $\Phi$ of weights, we study the following
classes: $\Phi^+$ of {\em $\Phi$-flat}
weights which are the $\psi$ for which $\psi$-colimits commute in the 
base $\V$ with limits with weight in $\Phi$; and 
$\Phi^-$ - dually defined - of weights 
$\psi$ for which $\psi$-limits commute in the 
base $\V$ with colimits with weight in $\Phi$.
We show that both these classes are {\em saturated} (i.e. 
{\em closed} with 
the terminology of \cite{AK88}).
We prove that for the class $\p$ of {\em all} weights
$\p^+$ = $\p^-$. 
For any small $\B$, we defined an enriched adjunction 
\`a la Isbell $[\B,\V]^{op} \rightharpoonup [{\B}^{op},\V]$ 
and show how it restricts to an equivalence 
${  ( \p^-({\B}^{op}) ) }^{op} \cong \p^-(\B)$ between
subcategories of {\em small projectives}. 
\end{abstract}

\begin{section}{Introduction}\label{section1}
The present notes had their beginning 
in the analysis of the results obtained by Borceux,
Quintero and Rosicky in their article
\cite{BQR98}. These authors were concerned with extending
to the enriched case the notion of accessibility
and its properties, described for ordinary categories
in the book \cite{MP89} of Makkai and Par\'e and 
\cite{AR94} of Adamek and Rosicky. They were led to discuss
categories - now meaning $\V$-categories - with finite 
limits (in a suitable sense), or more generally with 
$\alpha$-small ones, or with filtered colimits (in a suitable
sense), or more generally with $\alpha$-filtered colimits
and again to discuss the connexions between these classes
of categories. When we looked in detail at their work,
we observed that many of the properties they discussed
hold in fact for categories having colimits of {\em any}
given class $\Phi$, while other holds when $\Phi$ is the class
of colimits commuting in the base category $\V$ with the limits
of some class $\Psi$ - such properties as finiteness
or filteredness being irrelevant to the {\em general} results. 
Approaching in this abstract way, not generalisations of accessibility
as such,
 but the study of categories with colimits (or limits)
of some class, brings considerable simplifications.
A recent work by Adamek, Borceux, Lack and Rosicky
\cite{ABLR02} also treats a general notion of accessibility
parameterised by pair of families of {\em ordinary} limits 
and colimits that commute in $Set$. Rather than families
of ordinary -or conical limits- we consider {\em indexed
or weighted} limits and colimits. As a result of this generalisation,
we will avoid some technical complications even in the case when 
the base $\V$ is $Set$.
Although our positive results are quite few in number, 
their value may be judged by their extra light they cast on the 
results in \cite{BQR98} or \cite{ABLR02}.\\

We begin by reviewing and completing some known facts in 
the first sections \ref{section2}, \ref{section3}.
Section \ref{section4} contains a slight extension of Kelly's
cocompletion theorem.
Section \ref{section5} treats generally the limit/colimit commutation 
in the base $\V$ and contains a study of families of the form 
$\Phi^+$ of {\em $\Phi$-flat} weights - 
i.e. weights whose colimits commute with $\Phi$-weighted limits
and families of the form $\Phi^-$ of weights whose limits 
commute with $\Phi$-weighted colimits.
Section \ref{section6} focusses on the family $\p^-$ for $\p$
the class of all weights. $\p^-$ being the class of {\em small 
projective}, is shown to be also the class of $\p$-flat weights.
More results regarding the Cauchy-completion and the Isbell
adjunction in the enriched context are then presented.
\end{section}

\begin{section}{Revision of the terminology}
\label{section2} 
The necessary background knowledge about 
enriched categories is largely contained in \cite{Kel82}, 
augmented by Kelly's ``Amiens'' article \cite{Kel82-2}, and 
by the Albert-Kelly article \cite{AK88}.\\

We use ``category'', ``functor'', and ``natural transformation''
for ``$\V$-category'', ``$\V$-functor'', and ``$\V$-natural
transformation'', except when more precision is needed. As usual the ordinary
category $\V_o$ is supposed locally small, complete and cocomplete.
$\VCAT$ is the 2-category of $\V$-categories whereas $\VCat$ is
that of small $\V$-categories. $Set$ is the category of small
sets, $\Cat = \SCat$ is the 2-category of small categories,
$\CAT = \SCAT$ the 2-category of locally small categories.\\  

A {\em weight} is a functor $\phi: \K \rightarrow \V$ 
with domain $\K$ {\em small}; weights were called {\em indexing-types}
in \cite{Kel82} and in \cite{Kel82-2}. Recall that the 
{\em $\phi$-weighted
limit} $\{\phi, T\}$ of a functor $T:\K \rightarrow \A$ is defined
representably by
\begin{fact}\label{fac2.1}
$\A(a, \{ \phi, T \}) \cong [\K, \V](\phi, \A(a,T-))$,
\end{fact}
while the {\em $\phi$-weighted colimit} $\phi*S$ for 
$S:\K^{op} \rightarrow \A$ is defined dually by
\begin{fact}\label{fac2.2}
$\A(\phi * S,a) \cong [\K, \V](\phi, \A(S-,a))$,
\end{fact}
so that $\phi * S$ is equally the $\phi$-weighted 
limit of $S^{op}: \K \rightarrow \A^{op}$. Of course 
the limit $\{ \phi, T \}$ consists not just of the object
$\{ \phi, T \}$ but also of the representation \ref{fac2.1},
or equally of the corresponding {\em counit}
$\mu: \phi \rightarrow \A(\{\phi,T  \},T-)$; it is by 
{\em abus de langage}
that we usually mention only 
$\{\phi, T  \}$. When $\V = Set$, we refind the classical 
(or ``conical'') limit of $T:\K \rightarrow \A$ and the classical
colimit of $S:\K \rightarrow \A$ as 
\begin{fact}\label{fac2.3}
$lim\;T = \{ \Delta 1 , T \}$ $\;$and$\;$ $colim\;S = \Delta 1 * S$
\end{fact} 
where $\Delta 1: \K \rightarrow Set$ is the constant functor
at the one point set $1$.
Recall that the weighted limits and colimits can be calculated 
using the classical ones when $\V = Set$: for then the 
presheaf $\phi: \K \rightarrow Set$ gives the discrete
op-fibration $d: el(\phi) \rightarrow \K$ where
$el(\phi)$ is the category of elements of $\phi$, and now
\begin{fact}\label{fac2.4}
$\{ \phi, T \} = lim \{ \xymatrix{ el(\phi) \ar[r]^{d}&  \K \ar[r]^{T}
  & \A} \}$,
\end{fact}
\begin{fact}\label{fac2.5}
$ \phi * S = colim \{ \xymatrix{ {el(\phi)}^{op} \ar[r]^{d^{op}}&  
\K^{op} \ar[r]^{S} & \A}    \} $
\end{fact}
Also a functor $F: \B \rightarrow \C$ is said to {\em preserve
the limit} $\{ \phi, T \}$ as in \ref{fac2.1} when
$F(\{ \phi, T \})$ is the limit of $FT$ weighted by 
$\phi$ with counit 
$$\xymatrix{ \phi \ar[r]^(.3){\mu}
&  \A(\{ \phi, T \},T-) \ar[r]^(.45){F} 
& \B(F\{ \phi, T \},FT-) }$$ and $F$ is said to 
preserve the colimit
$\phi * S$ as in \ref{fac2.2} when $F^{op}$ preserves
$\{ \phi, S^{op} \}$.\\

We spoke above of a ``class $\Phi$ of colimits'' or a 
``class $\Psi$ of limits''; but this is loose and rather
dangerous language - the only thing that one can sensibly speak
of is {\em a class $\Phi$ of weights}. Then a category $\A$ {\em admits
$\Phi$-limits} or is {\em $\Phi$-complete}, if $\A$ admits the limit
$\{\phi, T \}$ for each weight $\phi: \K \rightarrow \V$ of
$\Phi$ and each $T: \K \rightarrow \A$, while $\A$ {\em admits 
$\Phi$-colimits}, or is {\em $\Phi$-cocomplete} when $\A$ admits 
the colimit $\phi * S$ for each $\phi: \K \rightarrow \V$
in $\Phi$ and for each $S: \K^{op} \rightarrow \A$
(and thus when $\A^{op}$ is $\Phi$-complete). Of course a 
functor $\A \rightarrow \B$ between $\Phi$-complete categories
is said to be {\em $\Phi$-continuous} when it preserves all 
$\Phi$-limits and one defines {\em $\Phi$-cocontinuous} dually.
We write $\PhiCont$ for the 2-category of $\Phi$-complete
categories, $\Phi$-continuous functors, and all natural
transformations - which is a (non full) sub-2-category
of $\VCAT$; and similarly $\PhiCo$ for the 2-category
of $\Phi$-cocomplete categories, $\Phi$-cocontinuous functors, 
and all natural transformations.\\
 
To give a class $\Phi$ of weights is to give for each small
$\K$, those $\phi \in \Phi$ with domain $\K^{op}$; let us use
as in \cite{AK88} the notation
\begin{fact}\label{fac2.6}
$\Phi[\K] = \{ \phi \in \Phi \mid dom (\phi) = \K^{op} \}$,
\end{fact}
so that
\begin{fact}\label{fac2.7}
$\Phi = \Sigma_{\K\;small} \Phi[\K]$.
\end{fact}
In future, we look at $\Phi[\K]$ as a full subcategory
of the functor category $[\K^{op},\V]$. The smallest class
of weights is the empty class $0$, and $\ZCont$ is just 
$\VCAT$. The largest class of weights consists of all
weights with small domains, and we denote this class
by $\p$; the 2-category $\PCont$ is just the 2-category
$\Cont$ of {\em complete} categories and {\em continuous} functors,
and similarly $\PCo = \Co$.\\ 

There may well be different classes $\Phi$ and $\Psi$
for which the sub-2-categories $\PhiCont$ and $\PsiCont$
of $\VCAT$ coincide; which is equally to say that 
$\PhiCo$ and $\PsiCo$ coincide.
When $\V = Set$, for instance, $\Cont = \PCont$ coincides
with $\PhiCont$ where $\Phi$ consists of the weights
for products and for equalisers. We define the 
{\em saturation} $\Phi^*$ of a class $\Phi$ of weights
as follows: the weight $\psi$ belongs to $\Phi^{*}$ when
every $\Phi$-complete category is also $\psi$-complete
and every $\Phi$-continuous functor is also 
$\psi$-continuous. 
Actually we change here the name given in \cite{AK88} where
$\Phi^*$ was called the {\em closure} of $\Phi$ but we wish to 
avoid overusing this terminology.
Clearly then, we have
\begin{fact}\label{fac2.8}
$\PhiCont = \PsiCont \Leftrightarrow \PhiCo = \PsiCo 
\Leftrightarrow \Phi^* = \Psi^*.$
\end{fact} 
When $\V = Set$, we can of course consider 
$\PhiCont$ where $\Phi$ consists of the $\Delta 1 : \K \rightarrow Set$
for all $\K$ in some small class $\D$ of small categories;
and we might write $\DCont$ for this 2-category $\PhiCont$
of {\em $\D$-complete} categories, {\em $\D$-continuous} functors, and 
all natural transformations. 
We underline the fact that however when $\V = Set$,
the $\PhiCont$ generally do NOT occur as $\DCont$ as 
above. A simple example of this situation is in \cite{AK88}.\\

There is a distinction which is often important, but which
is not covered by the notation above. We spoke of $\VCAT$
as a 2-category, the category $\VCAT(\A,\B)$ having as its 
objects the $\V$-functors $T:\A \rightarrow \B$ and as its 
arrows, the $\V$-natural transformations $\alpha: T \rightarrow S: \A
\rightarrow \B$. When $\A$ is small, however, we also have the 
$\V$-category $[\A,\B]$, whose underlying category $[\A,\B]_0$
is $\VCAT(\A,\B)$. Thus {\em small} $\V$-categories form not only 
a 2-category $\VCat$, but also a $\VCat$-category
(as $\VCat$ is a monoidal closed 2-category and as any
closed category, it is enriched over itself).
Again, even when $\A$ is not small, we can (as in
\cite{Kel82} section 3.12) see $[\A,\B]$ as a $\Vp$-category
where $\Vp$ is an extension of $\V$ to a $Set'$-category
(with $Set'$ a universe containing $Obj(Set)$ and $Obj(\B)$
for all categories $\B$), such that 
$\V \rightarrow \V'$ preserves all limits that exist in $\V$
and all $Set'$-small colimits that exist in $\V$.
So we take the view that $[\A,\B]$ always exists as 
a $\V'$-category, which is fact a $\V$-category when 
$\A$ is small. Thus, even for large $\A$, $\A$ is $Set'$-small
and we may see 
$\VCAT(\A,\B)$ as the ordinary category underlying 
the $\V'$-category $[\A,\B]$, so that $\V$-categories
form not only the 2-category $\VCAT$ cut also 
a $\VpCAT$-category denoted $\underline{\VCAT}$. 
This has a full subcategory
given by the $\Phi$-complete categories; and if $\A$ and $\B$
are $\Phi$-complete, the $\Phi$-continuous functors $\A \rightarrow \B$
determine a full subcategory $\PhiCont[\A,\B]$ of the 
$\V'$-category $[\A,\B]$. Note the use of {\em square} brackets in 
the name of the $\Vp$-category $\PhiCont[\A,\B]$ (which is a 
$\V$-category for small $\A$); its underlying category is our
earlier $\PhiCont(\A,\B)$ with round brackets.
Thus the 2-category $\PhiCont$ really underlies a 
$\VpCAT$-category that we might note $\underline{\PhiCont}$ with 
\begin{fact}\label{2.9}
$\underline{\PhiCont}(\A,\B) = \PhiCont[\A,\B]$
\end{fact} 
(Note that we suppress the mention of $\V$ in $\PhiCont$
or $\underline{\PhiCont}$; one cannot know the class $\Phi$ without knowing
$\V$. There may be a problem when $\Phi = \p$, so that 
$\PhiCont$ is just $\Cont$; one can write $\VCont$ if
necessary.)
Of Course the definition of $\underline{ \PhiCont }$ makes it clear 
that \ref{fac2.8} can equally be written as
\begin{fact}\label{2.10}
$\underline{ \PhiCont } = \underline{ \PsiCont }$ 
$\;\Leftrightarrow\;$ 
$\underline{ \PhiCo } = \underline{ \PsiCo }$ 
$\;\Rightarrow\;$ $\Phi^* = \Psi^*$.
\end{fact}
\end{section}

\begin{section}{Free cocompletions of categories and
saturated classes of weights}\label{section3} 
Another piece of background 
knowledge that we need to recall concerns  the ``left bi-adjoint'' to the 
forgetful 2-functor  $U_{\Phi}: \PhiCo \rightarrow \VCAT$,
or rather to 
$U_{\Phi}: \underline{\PhiCo} \rightarrow \underline{\VCAT}$.
(Note that it is convenient to deal with colimits
rather than limits.)\\

In the language of \cite{Kel82}, a presheaf $F:\A^{op} \rightarrow \V$
is said to be {\em accessible} if it is the left Kan extension
of its restriction to some small full subcategory of $\A^{op}$.
Even though $[\A^{op}, \V]$ is when $\A$ is not small, only
a $\V'$-category, its full subcategory $\p \A$ having as objects the accessible
presheaves is a $\V$-category, as follows from
\cite{Kel82} 4.41; of course $\p \A$ coincides with 
$[\A^{op},\V]$ when $\A$ is small. Every representable
$\A(-,a)$ is accessible, so that the Yoneda embedding
$Y:\A \rightarrow [\A^{op},\V]$ takes its values in $\p \A$.
Remember that the functor category $[\A^{op},\V]$ admits
all small colimits, these being formed pointwise in $\V$.
Each object of $\p \A$, being a left Kan extension of the 
given kind, is a small colimit of representables.
Recall from \cite{Kel82} p.154 
that given a full subcategory $\A$ of $\B$ and
a family $\Phi$ of weights, the {\em closure of $\A$ in $\B$
under $\Phi$-colimits} is the smallest full replete 
subcategory of $\B$ containing $\A$ and closed 
under the formation $\Phi$-colimits in $\B$.  
Since it is shown in \cite{Kel82} Proposition 5.34 that 
$\p \A$ is closed in $[\A^{op},\V]$ under small colimits,
we see that $\p \A$ is precisely the {\em closure of $\A$ in 
$[\A^{op},\V]$ under small colimits} (when $\A$ is embedded in 
$[\A^{op},\V]$ by the Yoneda embedding $Y$).
For any class $\Phi$ of weights, and any category $\A$,
write $\Phi(\A)$ for the closure of $\A$ in $\p A$
(or equally in $[\A^{op},\V]$) under $\Phi$-colimits,
with $Z:\A \rightarrow \Phi(\A)$ and 
$W: \Phi(\A) \rightarrow [\A^{op},\V]$ for the full inclusions,
so that $Y:\A \rightarrow [\A^{op},\V]$ is the composite 
$WZ$. We now reproduce (the main point of) 
\cite{Kel82} Theorem 5.35. The proof below is a 
little more direct than that given there, which referred
back to earlier results. The result itself must be older still, at
least for certain classes $\Phi$.

\begin{proposition}
\label{Pro3.1}
For any $\Phi$-cocomplete category $\B$ there is an
equivalence of $\Vp$-categories
\begin{fact}\label{fac3.1}
$[Z,1]:\PhiCo[\Phi(\A),\B] \rightarrow [\A,\B]$
\end{fact}
with an equivalence inverse given by $Lan_Z$.
It follows that $\PhiCo[\Phi(\A),\B]$ is like $[\A,\B]$,
a $\V$-category when $\A$ is small. Thus $\Phi(-)$ provides
a left bi-adjoint, not only to the forgetful 2-functor
$\PhiCo \rightarrow \VCAT$, but even  to the forgetful
$\VpCat$-functor $\underline{\PhiCo} \rightarrow \underline{\VCAT}.$
\end{proposition}
\pf
The left Kan extension 
$Lan_ZG$, if it exists, is given by 
$Lan_ZG(F) = \tilde{Z}F*G$ where
$\tilde{Z}(F) = \Phi(\A)(Z-,F)$. However
$\Phi(\A)(Z-,F) 
= [\A^{op}, \V](WZ-,WF) = [\A^{op},\V](Y-,WF)$,
which by Yoneda is isomorphic to $WF$. 
Consider the full subcategory of $\Phi(\A)$
given by those $F$ for which $WF*G$ does exist; it contains the 
representables by Yoneda, 
it is also closed under $\Phi$-colimits according to 
\cite{Kel82} (3.22)
since these colimits exist in $\B$; so it is all of $\Phi(A)$.
What is more: $Lan_Z(G) = W-*G: \Phi(\A) \rightarrow \B$
preserves $\Phi$-colimits as $W$ preserves $\Phi$-colimits
and $-*G$ is cocontinuous (straightforward from \cite{Kel82} (3.22) again).
So one does have indeed a functor 
$Lan_Z: [\A,\B] \rightarrow \PhiCo[\Phi(\A),\B]$, while one
has trivially the restriction functor 
$[Z,1]: \PhiCo[\Phi(\A),\B] \rightarrow [\A,\B]$.
The canonical $G \rightarrow Lan_Z(G)Z$ is invertible
for all $G$ since $Z$ is fully faithful by \cite{Kel82}
(4.22). So that it remains to consider the canonical
$\alpha: Lan_Z(SZ) \rightarrow S$ for a $\Phi$-cocontinuous 
$S: \Phi(\A) \rightarrow \B$. The $F$-component of
$\alpha$ for $F \in \Phi(\A)$ is the canonical 
$\alpha_F: WF * SZ \rightarrow SF$; and clearly
the collection of those $F$ for which $\alpha_F$ is 
invertible contains the representables and is closed 
under $\Phi$-colimits, therefore it is the totality
of $\Phi(\A)$.
\epf

\begin{remark}
\label{Rem3.2}
It is often said that $\Phi(\A)$ is the {\em free
$\Phi$-cocomplete category on $\A$}, but by this is 
usually meant an equivalence 
$\PhiCo(\Phi(\A),\B) \cong \VCAT(\A,\B)$ of categories
which is {\em weaker} than \ref{Pro3.1}.
\end{remark} 

\begin{remark}
\label{Rem3.3} As a particular case $\p \A$ itself is the free 
cocomplete category on $\A$; in other words $\Phi(\A) = \p \A$ 
when $\Phi$ is the class of {\em all} weights - which is why
(identifying $\p(\A)$ with $\p \A$) we use $\p$ as the name 
for this class of all weights.
\end{remark}

We can form $\Phi(\A)$ by transfinite induction. Defining 
successively full replete subcategories $\A_{\alpha}$ of $\p \A$
as $\alpha$ runs through the ordinals: $\A_0$ is the
category $\A \subseteq \p \A$ of representables, $\A_{\alpha+1}$
consists of $\A_{\alpha}$ together with all $\Phi$-colimits
in $\p \A$ of diagrams in $\A_{\alpha}$, and for a limit ordinal
we let $\A_{\alpha} = \bigcup_{\beta < \alpha} \A_{\beta}$,
then $\Phi(\A)$ is the smallest full replete subcategory of 
$[A^{op},\V]$ containing all the $\A_{\alpha}$ for all ordinals 
$\alpha$.  
Note that this shows that {\em $\Phi(\A)$ is small when $\A$  
and $\Phi$ are small}. However $\Phi(\A)$ may be small
for all small $\A$ when $\Phi$ is not a small class: in this 
case we call the class $\Phi$ {\em locally small} (as usual in category
theory, we call a category {\em small} when its set of 
{\em isomorphism classes} is a small set).
In a number of important cases, one has $\Phi(\A) = \A_1$ in the 
notation above; but there is no special value in this condition,
which (as we shall see) always holds for a small $\A$ when 
the class $\Phi$ is {\em saturated}.\\ 

An explicit description of the saturation $\Phi^*$ of a class
$\Phi$ of weights was given by Albert and Kelly in \cite{AK88},
in the following terms:
\begin{proposition}
\label{Pro3.4} The weight $\psi: \K^{op} \rightarrow \V$ lies in the
saturation $\Phi^*$ of the class $\Phi$ if and only if the 
object $\psi$ of $\p \K = [\K^{op},\V]$ lies in the full 
subcategory $\Phi(\K)$ of $[\K^{op}, \V]$.
\end{proposition}
There is another useful way of putting this. When $\K$ is small,
both $\Phi[\K]$ and $\Phi(\K)$ makes sense for any class 
$\Phi$; and in fact we have
\begin{fact}
\label{F3.2} $\Phi[\K] \subseteq \Phi(\K),$
\end{fact}
since for $\phi: \K^{op} \rightarrow \V$ the Yoneda 
isomorphism 
\begin{fact}
\label{F3.3}
$\phi \cong \phi * Y$
\end{fact} exhibits $\phi$ as 
an object of $\Phi(\K)$ when $\phi \in \Phi$.
We can write Proposition \ref{Pro3.4} as
\begin{fact} \label{F3.4}
$\Phi^*[\K] = \Phi(\K)$
\end{fact}
so that $\Phi$ is a saturated class precisely 
when 
\begin{fact} \label{F3.5}
$\Phi[\K] = \Phi(\K)$
\end{fact}
for each small $\K$. In other words the class $\Phi$ is 
saturated precisely when, for each small $\K$, {\em the full subcategory
$\Phi[\K]$ of $[\K^{op},\V]$ contains the representables 
$\K(-,k)$ and is closed in $[\K^{op}, \V]$ under $\Phi$-colimits.}

Note that the saturation $0^*$ of the empty class $0$ consists
precisely in the representables - now in the sense that $0^*[\K] =
0^*(\K)$ 
consists exactly in the functors {\em isomorphic 
to} some $\K(-,k): \K^{op} \rightarrow \V$. 

It follows of course from the definitions of $\Phi(\A)$
and of $\Phi^*$ that 
\begin{fact}\label{F3.6}
$\Phi^*(\A) = \Phi(\A)$ 
\end{fact}
for any $\A$, small or not. We cannot write \ref{F3.5}
when $\K$ is replaced by a non-small $\A$, since then $\Phi[\A]$
has no meanings; but a partial replacement for it is provided
by the following which was Proposition 7.4 in \cite{AK88}
\begin{proposition}
If the presheaf $F:\A^{op} \rightarrow \V$ lies in $\Phi(\A)$ for 
some closed class $\Phi$ then $F$ is a $\Phi$-colimit of 
representables; that is $F \cong \phi * ZS$ for some
$\phi: \K^{op} \rightarrow \V$ in $\Phi$ and some
$S: \K \rightarrow \A$.   
\end{proposition}

It may be useful to understand extreme special cases
of one's notation: we mentioned above the empty class
$0$ of weights, with $0[\K] = 0$ and 
with $0(\K) = 0^{*}(\K) = 0^*[\K]$ consisting of the representables.
Another extreme case involves the empty $\V$-category $0$
with no objects. Of course $P0 = [0^{op},\V]$ is the 
terminal category $1$; its unique object is the unique 
functor $!:0^{op} \rightarrow \V$ and $1(!,!)$ is 
the terminal object $1$ of $\V$. (This differs in general
from the {\em unit $\V$-category} $\I$, with one object $*$
but with $\I(*,*)$$=$$I$, the unit for $\otimes$.) 
So for any class $\Phi$, we have $\Phi[0] = 0$ if 
$!: 0^{op} \rightarrow \V$ is not in $\Phi$, and 
$\Phi[0] = 1$ otherwise. Now $\Phi(0)$ is the closure
of $0$ in $P0$ under $\Phi$-colimits, and any diagram
$T: \K \rightarrow 0$ has $\K=0$, so that $\Phi(0) = 0$
if $! \not \in \Phi$ and otherwise $\Phi(0)$ contains
$ ! * Y = !$, giving $\Phi(0) = 1$. So in fact 
$\Phi(0) = \Phi[0]$, being $0$ or $1$. Both are possible for a 
closed $\Phi$ for $P0=1$, while the Albert-Kelly theorem
(Proposition \ref{Pro3.4})  gives $0^*[0] = 0(0)= 0[0] = 0$.\\

Before ending this section, we recall a result characterising 
$\Phi$-cocomplete categories, along with a short proof. This 
was Proposition 4.5 in \cite{AK88}.
\begin{proposition}
For any class $\Phi$ of weights, a category $\A$ admits 
$\Phi$-colimits if and only if the fully faithful
embedding $Z:\A \rightarrow \Phi(\A)$ admits a left adjoint;
that is, if and only if the full subcategory $\A$ given by
the representables is reflective in $\Phi(\A)$.
\end{proposition}
\pf
If $\A$ is reflective, it admits $\Phi$-colimits
because $\Phi(\A)$ does so. Suppose conversely that 
$\A$ admits $\Phi$-colimits, and write $\B$ for the full
subcategory of $[\A^{op},\V]$ given by those objects
admitting a reflection into $\A$; then $\B$ contains $\A$
and $\B$ is closed in $[\A^{op},\V]$ under $\Phi$-colimits
since $\A$ admits these; so that $\B$ contains $\Phi(\A)$,
as desired.
\epf
\end{section}

\begin{section}{The embedding $\A \rightarrow \Phi(\A)$}\label{section4} 
We recall from Proposition 5.62 of \cite{Kel82}
a result characterising categories of the form $\Phi(\A)$
- or more precisely functors of the form 
$Z: \A \rightarrow \Phi(\A)$ (which we also write when necessary
as $Z_{\A}$). At the same time, we give a direct proof; for the proof
in \cite{Kel82} refers back to earlier results in that book.

We begin with a piece of notation: for a category $\A$ and a 
class $\Phi$ of weights, we write $\A_{\Phi}$ for the full 
subcategory of $\A$ generated by those $a \in \A$ for which the
representables $\A(a,-): \A^{op} \rightarrow \V$ preserve all
$\Phi$-colimits (That is all $\Phi$-colimits that {\em exist in}
$\A$). There is no agreed name for $\A_{\Phi}$; the objects of
$\A_{\Phi}$ are usually called {\em finitely presentable}
when the $\Phi$-colimits are the classical filtered colimits;
while when $\Phi$ is the class $P$ of all weights the objects 
of $\A_\Phi$ were called {\em small projective} in \cite{Kel82},
but have also been called {\em atoms} by some authors. 
Therefore we propose to name the objects of $\A_{\Phi}$ 
{\em $\Phi$-atoms}.
When $\A$ admits $\Phi$-colimits and hence $\Phi^*$-colimits,
it follows from the definition of $\A_{\Phi}$ that
\begin{fact}\label{fac4.1}
$\A_{\Phi^*} = \A_{\Phi}.\label{F4.1}$ 
\end{fact}
The following is the characterisation result of Proposition
$5.62$ \cite{Kel82} with a slightly expanded form of its 
statement.
\begin{proposition}\label{Pro4.1}
In order that $G:\A \rightarrow \B$ be equivalent to the free
$\Phi$-cocompletion $Z:\A \rightarrow \Phi(\A)$ of $\A$ for a
class $\Phi$ of weights, the following conditions are 
necessary and sufficient:
\begin{itemize}
\item $(i)$ G is fully faithful (allowing us to treat $\A$
henceforth as a full subcategory of $\B$); 
\item $(ii)$ $\B$ is $\Phi$-cocomplete;
\item $(iii)$ the closure of $\A$ in $\B$ under $\Phi$-colimits
is $\B$ itself;
\item $(iv)$ $\A$ is contained in the full subcategory $\B_{\Phi}$
of $\Phi$-atoms of $\B$.
\end{itemize}
When these conditions hold, the functor 
$\tilde{G}: \B \rightarrow [\A^{op},\V]$ defined (
in the usual notation of \cite{Kel82}) by 
$\tilde{G}b = \B(G-,b)$ is fully faithful and takes its image
in the full subcategory $\Phi(\A)$ of $[\A^{op},\V]$ thereby defining
a fully faithful $K: \B \rightarrow \Phi(\A)$; and this $K$ is 
in fact an equivalence - inverse being given by
$Lan_Z(G): \Phi(\A) \rightarrow \B$, which by Proposition \ref{Pro3.1}
is the unique $\Phi$-cocontinuous extension of $G$ to 
$\Phi(\A)$.
\end{proposition}
\pf The necessity of the first three conditions is clear;
for that of the forth, writing as before
$$\xymatrix{\A \ar[rd]_{Z} \ar[rr]^{Y} & &  [\A^{op},\V] \\
& \Phi(\A) \ar[ru]_{W}  &}$$ for the full inclusions, the
point is that $\Phi(\A)(Za,-) \cong [\A^{op},\V](Ya,W-)$ preserves
$\Phi$-colimits since $W: \Phi(\A) \rightarrow [\A^{op},\V]$
does so by definition of $\Phi(\A)$, while 
$[\A^{op},\V](Ya,-): [\A^{op},\V] \rightarrow \V$ preserves
all small colimits, being isomorphic by Yoneda to the 
evaluation $E_a$.

Turning to the proof of sufficiency, to say that $\tilde{G}$
is fully faithful is just to say (see Theorem 5.1 \cite{Kel82})
that $G: \A \rightarrow \B$ is dense; and in the present case
that is so by Theorem 5.19 of \cite{Kel82} since $\B$ is the closure
of $\A$ in $\B$ under $\Phi$-colimits and $\A$ lies in $\B_{\Phi}$.
The proof is so short that we may recall it: we consider
the full subcategory of $\B$ given by those $b$ for which
$\tilde{G}_{b,c}: \B(b,c) \rightarrow [\A^{op},\V](\tilde{G}b,\tilde{G}c)$ 
is invertible for all $c$, observing that it is all of $B$ since 
it contains $\A$ - since $G$ is fully faithful - and is closed 
under $\Phi$-colimits - as $\A \subseteq \B_{\Phi}$.

Next, to see that each $\tilde{G}b = \B(G-,b)$ lies in the 
full replete subcategory $\Phi(\A)$ of $[\A^{op},\V]$, consider
the full subcategory of $\B$ given by those $b$ for which
this is so; this contains $\A$ since $\tilde{G} \cdot Ga \cong Ya$
because $G$ is fully faithful, and it is closed in $\B$ under 
$\Phi$-colimits since $\tilde{G}$ preserves these by 
$(iv)$ and since $\Phi(\A)$ is closed under these in 
$[\A^{op},\V]$; so it is all of $\B$. Thus $\tilde{G}$ is indeed 
of the form $\tilde{G} = WK$ for some fully faithful
$K: \B \rightarrow \Phi(\A)$.

It remains to show that $K$ and 
$S = Lan_{Z}G: \Phi(\A) \rightarrow \B$ are equivalence-inverses.
Recall from Proposition \ref{Pro3.1} that $S$ is equally the
essentially unique $\Phi$-cocontinuous functor with $SZ = G$. Since 
$\tilde{G} = WK$ and $G \cong SZ$, one has $WKSZ \cong \tilde{G}G 
\cong Y \cong WZ$, giving $KSZ \cong Z$ since $W$ is fully faithful
and then giving $KS \cong 1$ by Proposition $\ref{Pro3.1}$ since 
$KS$ and $1$ are $\Phi$-cocontinuous (K being so because 
$WK = \tilde{G}$ is so).
Finally $KS \cong 1$ gives $KSK \cong K$ and $SK \cong 1$
since $K$ (as we saw) is fully faithful.
\epf
This is of course of particular interest in the case of a small
$\A$ in Proposition \ref{Pro4.1}. If we suppose the class $\Phi$ 
saturated for simplicity, we may cast the result for small $\A$'s 
in the form 
\begin{proposition}\label{Pro4.2}
For a saturated class $\Phi$ of weights, the following properties
of a category $\B$ are equivalent:
\begin{itemize}
\item $(i)$ For some small $\K$, there is an equivalence $\B \cong \Phi(\K)$;
\item $(ii)$ $\B$ is $\Phi$-cocomplete and has a small full subcategory
$A \subseteq \B_{\Phi}$ such that every object of $\B$ is a 
$\Phi$-colimit of a diagram in $\A$;
\item $(iii)$ $\B$ is $\Phi$-cocomplete and has a small full
subcategory $\A \subseteq \B_{\Phi}$ such that the closure
of $\A$ in $\B$ under $\Phi$-colimits is $\B$ itself.
\end{itemize}
Under the hypothesis $(iii)$ - and so a fortiori under $(ii)$ -
if $G: \A \rightarrow \B$ denotes the inclusion, the functor 
$\tilde{G}: \B \rightarrow [\A^{op},\V]$ is fully faithful,
with $\Phi(\A)$ for its replete image.
\end{proposition}

\begin{remarks}\label{Rem4.3}
\begin{itemize}
\item (a) When $\Phi$ is the class $\p$ of all weights, we get a
characterisation here of the functor category $\p \K = [\K^{op},\V]$
for a small $\K$; note that it differs from the characterisation given
in \cite{Kel82} Theorem 5.26 which replaces the condition that 
$\B$ be the colimit closure of $\A$ by the condition that $\A$ be 
strongly generating in $\B$; but that these conditions are very similar
in strength by \cite{Kel82} Proposition 3.40.
\item (b) When $\V= Set$, various special cases of this are well known,
besides that where $\Phi = \p$; for instance the cases where 
$\Phi$-colimits are finite coproducts, or finite colimits 
or filtered colimits, or absolute colimits. We will return to 
this last example in section \ref{section6}.
\item (c) Theorem 2.4 of \cite{BQR98} is the special case where $\V$ is 
locally finitely presentable as a closed category in the sense
of \cite{Kel82-2} and $\Phi$ is the saturated class of flat 
presheaves.  
\end{itemize}
\end{remarks}
$\;$\\

Of particular interest regarding Proposition \ref{Pro4.2} are the saturated
families of the form $\Psi^{+}$ and $\Phi^-$, for arbitrary classes
of weights $\Psi$ and $\Phi$.  
In particular we shall see that when $\Phi = \p$ is the family 
of all weights, $\Phi^+(\A)$ is the Cauchy-completion 
of $\A$.
Also the fact of considering families of weighted limits/colimits
instead of mere conical limits/colimits is crucial even in the case
$\V=Set$. In \cite{ABLR02}, Adamek, Borceux, Lack and Rosicky
investigated the notion of general commutation of
families of limits/colimits - the limits and colimits
considered are classical - i.e.defined in terms of diagrams 
- not weight. To obtain a result similar to \ref{Pro4.2} (that one may 
call ``parameterised  accessibility''), they needed to impose the
technical
condition of ``soundness'' on classes of limits. The next
section should be a key argument in favour of weighted limits/colimits, that
should also invalidate the wrong idea that ``weighted limits/colimits
in Sets are just the same as usual ones''.\\
\end{section}

\begin{section}{Limits and colimits commuting in $\V$}\label{section5}
The new observations to which we now turn on begin with the general
theory of the commutativity in $\V$ of limits and colimits.
For a pair of weights $\psi: \K^{op} \rightarrow \V$ and 
$\phi: \LL^{op} \rightarrow \V$, to say that
\begin{fact}\label{fac5.1} 
$\phi * - : [\LL,\V] \rightarrow \V$ preserves $\psi$-limits
\end{fact}
is equally to say that 
\begin{fact}\label{fac5.2} 
$\{ \psi, - \} : [\K^{op},\V] \rightarrow \V$ preserves $\phi$-colimits
\end{fact}
because each in fact asserts the invertibility, for every functor 
$S: \K^{op} \otimes \LL \rightarrow \V$, of the canonical 
comparison morphism
\begin{fact}\label{fac5.3}
$\phi? * \{ \psi-, S(-,?) \} \rightarrow \{ \psi-, \phi * S(-,?) \}.$
\end{fact}
$\;$\\
The reader who would not find
the above statement immediate, would be happy to consider 
the following proposition. 

\begin{proposition}
\label{commut}
Let $F: \A^{op} \rightarrow \V$,
$P: \K \rightarrow \V$ and $G: \A \otimes \K \rightarrow \V$
equivalent to $G': \A \rightarrow [\K,\V]$ and also 
to $G'': \K \rightarrow [\A,\V]$. Then   
$F*-: [\A, \V] \rightarrow \V$ preserves the limit 
$\{P , G''\}$ if and only if 
$\{P,-\} : [\K, \V] \rightarrow \V$ preserves the
colimit $F* G'$.
\end{proposition}
\pf 
Let $$\xymatrix{ P \ar[r]^-{\eta} & [\A, \V](\{ P, G''\},G'' - )   }$$
be the unit of $\{ P, G'' \}$ and 
$$\xymatrix{ F \ar[r]^-{\lambda} & [\K, \V](G'-, F*G') }$$
be the unit of $F*G'$.
We need to show that 
$$(1)\;\;\xymatrix{ 
P \ar[r]^-{\eta} & 
[\A, \V](\{ P, G''\},G'' -) \ar[r]^-{F*-} &
[F* \{P, G''\}, (F*G''-)] }$$
exhibits $F* \{P, G''\}$ as $\{P, F*G' \}$ if and only if 
$$(2)\;\;\xymatrix{ 
F \ar[r]^-{\lambda} & 
[\K, \V](F* G',G' -) \ar[r]^-{\{ P,- \}} & 
[\{P, F*G'\}, \{P, G'- \}]  }$$ exhibits
$\{P, F*G' \}$ as $F* \{P, G''\}$.\\

First note that given $x \in \V$, any natural in $k$, 
$$(1')\;Pk \rightarrow_k [x, (F*G')k]$$
corresponds via Yoneda to a natural in $v$,
$[v,x] \rightarrow_{v} [\K,\V](P,[v,(F*G')-])$.
Since $$[\K,\V](P,[v,(F*G')-] \cong_v [v,\{P, F*G' \}],$$
it corresponds also to an arrow $$(1'')\;\;x \rightarrow \{P, F*G' \}.$$
Also that $(1')$ exhibits
$x$ as the limit $\{P, F*G' \}$ is equivalent to 
the fact that $(1'')$ is iso.
Analogously any natural in $a$, $$(2')\;\; Fa \rightarrow [\{ P, G''\}a, x]$$ 
corresponds to an arrow $$(2'')\;\;F * \{P, G''  \} \rightarrow x,$$ and
$(2')$ exhibits $x$ as the colimit $F* \{P, G''  \}$ if and only if
$(2'')$ is iso.\\

Now the result follows from the fact the arrow $(1)$ above
corresponds by the bijection $(1')-(1'')$ 
to the same arrow as $(2)$ by $(2')-(2'')$.
This is now left to the reader.\\

When these statements are true for every such $S$, we say
that {\em $\phi$-colimits commute with $\psi$-limits in $\V$}.
For classes $\Phi$ and $\Psi$ of weights, if \ref{fac5.1} (equivalently
\ref{fac5.2}) holds for all $\phi \in \Phi$ and all $\psi \in \Psi$,
we say that {\em $\Phi$-colimits commute with $\Psi$-limits in $\V$}.
For any class $\Psi$ of weights we may 
consider the class $\Psi^+$ of all weights $\phi$ for which 
$\phi$-colimits commute with $\Psi$-limits in $\V$; and for any class
$\Phi$ of weights we may consider the class $\Phi^-$ of all weights
$\psi$ for which $\Phi$-colimits commute with $\psi$-limits
in $\V$. This forms of course a Galois connection 
(with $\Phi \subseteq \Psi^+$ if and only if 
$\Psi \subseteq \Phi^-$).
%
Note that $[\LL, \V]$ and $\V$ in \ref{fac5.1} admit all 
(small) limits; so that by the definition above 
of the saturation $\Psi^*$ of a class $\Psi$ of weights
if $\phi * -$ preserves all $\Psi$-limits, it also 
preserves $\Psi^{*}$-limits. From this and a dual 
argument, one concludes that:
\begin{proposition}\label{pro5.1}For any class
$\Phi$ and $\Psi$ of weights, the class $\Phi^-$ and 
$\Psi^+$ are saturated classes of weights; so that 
$\Psi^{+\;*} = \Psi^+$ and $\Phi^{-\;*} = \Phi^-$.
Moreover $\Psi^+ = \Psi^{*\;+}$ and
$\Phi^- = \Phi^{*\;-}$.
\end{proposition} 

When $\Psi$ consists of the weights for finite limits
(in the usual case for ordinary categories, or in the sense
of \cite{Kel82-2} when $\V$ is ``locally finitely presentable 
as a closed category'', its has been customary to call the 
elements of $\Psi^+$ the {\em flat} weights, as they
are those $\phi$ having $\phi * -: [\LL,\V] \rightarrow \V$
left exact. We might therefore for a general $\Psi$ think 
of the elements of $\Psi^+$ as {\em $\Psi$-flat} weights.   

Recall that the limit functor $\{\psi, - \}$ of 
\ref{fac5.2} is just the representable functor
$[\K^{op},\V](\psi,-):[\K^{op},\V] \rightarrow V$.
Accordingly $\psi:\K^{op} \rightarrow \V$ lies
in ${\Phi}^-$ for a given class $\Phi$ if and 
only if it lies in the subcategory $[\K^{op},\V]_{\Phi}$
of $\Phi$-atoms:
\begin{fact}\label{fac5.6}
$\Phi^-[\K] = \Phi^-(\K) = [\K^{op},\V]_{\Phi}$ 
\end{fact}
Part of the saturation of ${\Phi}^-$ - namely the 
closedness of ${\Phi}^-[\K]$ in $[\K^{op}, \V]$ under
$\Phi^-$-colimits is the special case for
$\A = [\K^{op},\V]$ of the following more general result:
\begin{proposition}\label{pro5.2} For any class $\Phi$ of weights 
and any category $\A$, the full subcategory $\A_{\Phi}$
of $\A$ is closed in $\A$ under $\Phi^-$-colimits 
that exists in $\A$.
\end{proposition}
\pf Let the colimit $\psi * S$ exist, where 
$\psi:\K^{op} \rightarrow \V$ lies in $\Phi^-$ and
$S:\K \rightarrow \A$ takes its values in $\A_{\Phi}$.
Then by definition 
$$\A(\psi * S,a) \cong [\K^{op}, \V](\psi, \A(S-,a)).$$
Since each $\A(SK,-)$ preserves $\Phi$-colimits,
and since $[\K^{op},\V](\psi,-)$ preserves $\Phi$-colimits
by \ref{fac5.6}, it follows that $\A(\psi*S,-)$ preserves
$\Phi$-colimits, that is to say $\psi*S \in \A_{\Phi}$. 
\epf

Let us add ad-hoc proof
of the closeness of $\Phi[\K]$ under $\Phi^+$-colimits in
$[\K^{op},\V]$ (though this does 
not bring anything more than Proposition \ref{pro5.1}) 
\begin{lemma}
For any $\Psi$-flat $G: \A^{op} \rightarrow \V$ and any functor
$H: \A \rightarrow [\C^{op},\V]$ with values $\Psi$-flat functors,
the colimit $G*H$ is again  $\Psi$-flat.       
\end{lemma}
\pf
Consider $F:\K \rightarrow \V$ in $\Psi$ and 
$L: \K \rightarrow [\C, \V]$. One has the successive
isomorphisms:
\begin{tabbing}
$(G * H) * \{ F,L \}$
\=$\cong$ 
\=$G * (H-* \{ F,L \})$ (
\cite{Kel82}, (3.23) ``continuity of a colimit in its index'')\\
\>$\cong$ 
\>$G * (\{ F, H- * L- \})$ 
\=(since for all $a$, $Ha*-$ preserves $\Psi$-limits)\\
\>$\cong$ 
\>$\{ F, G * (H- * L-) \}$
(since $G*-$ preserves  $\Psi$-limits)\\
\>$\cong$ 
\>$\{ F, (G * H)*L- \})$ (\cite{Kel82}, (3.23)).
\end{tabbing}
The resulting isomorphism 
$(G * H) * \{ F,L \}  \cong \{ F, (G * H)*L-  \})$
corresponds actually to the preservation of $\{F, L \}$ by
$(G*H)*-$.
\end{section}

\begin{section}{Cauchy completion and Isbell adjunction}\label{section6}
This last section is devoted to the study of the particular 
class of weights $\Q = \p^-$ where
$\p$ stands as before for the class of all weights.
That is to say that $\Q$ is the class defined by:
\begin{fact}
$\phi \in \Q$ $\Leftrightarrow$ $\phi$-limits commute with 
{\em all small} colimits in $\V$.  
\end{fact}
and equivalently, as we shall see, that the weights in 
$\Q$ are the so-called {\em small projective} ones.
According to Proposition \ref{pro5.1}, we already know that $\Q$ is a 
saturated class. We shall establish  
the following alternative characterisations for $\Q$. 
From Proposition \ref{sladj2}:
\begin{fact}
$\phi: \A^{op} \rightarrow \V$ lies in $\Q$ if and only 
the corresponding module $\xymatrix{ \I \ar[r]|{\circ} & \A}$
is left adjoint.
\end{fact}
This last result shows that the Cauchy-completion
of a small category $\A$ is just $\Q(\A) = \Q[\A]$.
Also from Proposition \ref{yac}:
\begin{fact} 
$\Q$ is the class of $\p$-flat weights, i.e. $\Q = \p^+$.
\end{fact}
The end of the section treats an enriched adjunction 
\`a la Isbell that restricts to an 
equivalence between subcategories of small projectives.\\ 

We need to recall a little terminology and few facts regarding
adjoints, Kan extensions and liftings in a bicategory.
Let us consider a bicategory $\bicatB$ with
horizontal composition ``$*$'' and
vertical one ``$.$'' to fix the notations.
$\bicatB^{op}$ 
will denote the bicategory obtained from $\bicatB$ by reversing 
arrows, and $\bicatB^{coop}$ will be the one obtained 
by reversing arrows and 2-cells.
We borrow from \cite{StWa78} the following terminology. 
Given a diagram as below in $\bicatB$
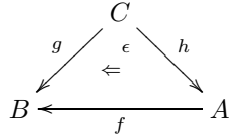
\begin{center}
$\xymatrix{
&   C \ar[ld]_g \ar[rd]^h &  \\
 B & &   A \ar[ll]^f \ar@{}[llu]^{\Leftarrow}_{\epsilon}
}$
\end{center}
the arrow $h$ together with 2-cell $\epsilon$ constitutes
the {\em right lifting of $g$ through $f$} if and only if 
pasting with $\epsilon$ at $h$ determines a bijection of 
2-cells
\begin{center} 
$\bicatB(C,B)(f*k,g) \longleftrightarrow \hskip -10pt \mid
\hskip +10pt \bicatB(C,A)(k,h).$
\end{center} 
In this situation we shall loosely write $h = \Rl(f,g)$.
The right lifting is said {\em respected} by a $k:X \rightarrow  C$
when the 2-cell $\epsilon * k$ exhibits $h * k$ as the 
right lifting of $g * k$ through $f$, it is {\em absolute}
when it is respected by any such arrow.
Note that a right lifting in $\bicatB$
is a right Kan extension in $\bicatB^{op}$ or a left Kan extension
in $\bicatB^{coop}$. One has also 
a dual notion of {\em left liftings} that correspond
to left Kan extensions in $\bicatB^{op}$.\\

The following result characterising
adjoints in a bicategory goes back to early work
of Lawvere and B\'enabou on closed bicategories. In fact 
the full result is really present in the one object-closed-bicategory 
situation given by a monoidal closed category, 
and that was studied in the symmetric case by \cite{KeLa80}.

\begin{proposition}\label{sladj1} In ANY bicategory, the following 
statements are equivalent:
\begin{itemize}
\item $(i)$ the arrow $s: A \rightarrow B$ has a right adjoint;
\item $(ii)$ for all $g: C \rightarrow B$, the right lifting 
$\Rl(s,g)$ of $g$ through $s$ exists and is absolute;
\item $(iii)$ the right lifting $\Rl(s,1)$ of 
$1:B \rightarrow B$ through $s$ exists and is respected by $s$.
\end{itemize}
When $s$ has a right adjoint, this one is $\Rl(s,1)$ the right lifting
of $1$ through $s$.
\end{proposition}
\pf $(i) \Rightarrow (ii)$. Let $\bicatB$ denote the ambient
bicategory. Let also 
$\eta: 1 \Rightarrow t*s: A \rightarrow A$
and 
$\epsilon: s * t \Rightarrow 1: B\rightarrow B$
denote respectively the unit and counit of $s \dashv t$.
For any 1-cells $g: C \rightarrow B$ and $x: C \rightarrow A$
the ``pasting with $\epsilon$'':
$( \epsilon*g) \cdot (s*-) : \bicatB(C,A)(x,t*g) \rightarrow
\bicatB(C,B)(s*x,g)$
admits for inverse the ``pasting with $\eta$'':
$(t*-) \cdot ( \eta * x )$. This shows that
$\epsilon * g: s * t * g \Rightarrow g$ exhibits $t*g$ as the right lifting
of $g$ through $s$. This lifting is trivially absolute.\\
$(ii)$ trivially implies $(iii)$.\\
$(iii) \Rightarrow (i)$.\\
The 1-cell 
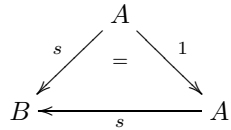
\begin{center}
$\xymatrix{& A \ar[ld]_s \ar[rd]^1 \ar@{}[d]|= & \\
B & & A \ar[ll]^s}$
\end{center}
factorises as
\begin{center}
$
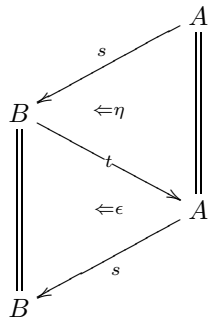
\xymatrix{ & & A \ar[lld]_s \ar@{=}[dd]\\
B \ar[rrd]|t \ar@{=}[dd] \ar@{}[rr]|{\Leftarrow \eta} & & \\
\ar@{}[rr]|{\Leftarrow \epsilon} & & A \ar[lld]^s\\
B & & }$
\end{center}
since $s$ respects the absolute right lifting
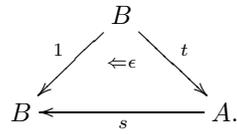
\begin{center}
$\xymatrix{& B \ar[ld]_1 \ar[rd]^t \ar@{}[d]|{\Leftarrow \epsilon}& \\
B & & A. \ar[ll]^s}$
\end{center}     
Now 
\begin{center}
$
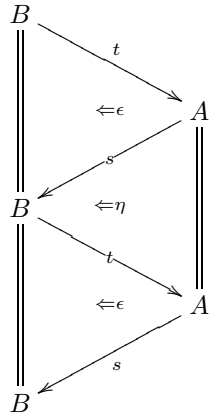
\xymatrix{B \ar@{=}[dd] \ar[rrd]^t && \\
\ar@{}[rr]|{\Leftarrow \epsilon}&& A \ar[lld]|s \ar@{=}[dd]\\
 B  \ar[rrd]|t  \ar@{}[rr]|{\Leftarrow \eta} \ar@{=}[dd] && \\
\ar@{}[rr]|{\Leftarrow \epsilon}& & \ar[lld]^s A \\
B &&
}$
\end{center}
equals 
\begin{center}
$
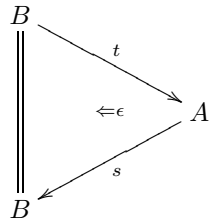
\xymatrix{B \ar@{=}[dd] \ar[rrd]^t&  &\\
\ar@{}[rr]|{\Leftarrow \epsilon} & & A \ar[lld]^s \\
B & & }$
\end{center}
and because $\epsilon$ is a right lifting
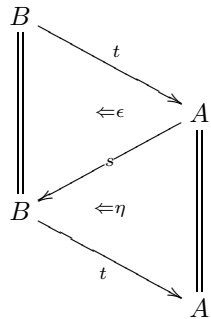
\begin{center}
$\xymatrix{
B \ar@{=}[dd] \ar[rrd]^t & & \\
\ar@{}[rr]|{\Leftarrow \epsilon} & & A \ar[lld]|s \ar@{=}[dd]\\
\ar@{}[rr]|{\Leftarrow \eta} B \ar[rrd]_t & & \\
& &  A}$
\end{center}
equals
\begin{center}
$
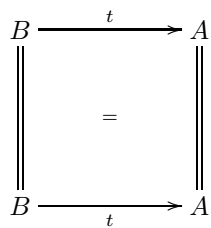
\xymatrix{B \ar@{=}[dd] \ar[rr]^t & &  A  \ar@{=}[dd] \\
\ar@{}[rr]|{=} & & \\
B  \ar[rr]_t & & A 
}$
\end{center}
As adjoints $F \dashv G$ in a bicategory $\bicatB$
correspond to adjoints $G \dashv F$ in $\bicatB^{op}$,
there are dual statements obtained by using left or right 
Kan extensions or left liftings. In detail
\begin{proposition} In any bicategory, the following assertions
are equivalent:
\begin{itemize}\label{sladj2} 
\item the arrow $s: A \rightarrow B$ has a right adjoint;
\item for all $g: A \rightarrow C$, the left Kan extension
$lan_s(g)$ exists and is absolute;
\item the left Kan extension  $lan_s(1)$ of $1: A \rightarrow A$ 
along $s$ exists and is preserved by $s$.
\end{itemize}
When $s$ has a right adjoint, this one is $lan_s(1)$.
\end{proposition}
And
\begin{proposition}\label{tradj} 
In any bicategory, the following assertions
are equivalent:
\begin{itemize} 
\item $(i)$ the arrow $t: B \rightarrow A$ has a left adjoint;
\item $(ii)$ for all $f: C \rightarrow A$, the left 
lifting $\Ll(t , f)$ of $f$ through $t$ exists and is absolute;
\item $(iii)$ the lifting $\Ll(t,1)$ of $1: A \rightarrow A$ 
through $t$ exists and is respected by $t$.
\item $(iv)$ for all $g: B \rightarrow C$, the right 
Kan extension $ran_t(g)$ of $g$ along $t$ exists and is 
absolute;
\item $(v)$ the right Kan extension $ran_t(1)$ of $1:B \rightarrow B$ 
along $t$ exists and is preserved by $t$.
\end{itemize}
When $t$ has a left adjoint, this one is $\Ll(t,1)$ or also
$ran_t(1)$. 
\end{proposition}
There is one lemma that we shall use later
\begin{proposition}
\label{lemadj}
In a bicategory if an arrow  $f: A \rightarrow B$
has a right adjoint $g: B \rightarrow A$
then for any $h:C \rightarrow B$, the right 
lifting of $h$ through $f$ exists and is given by
$\Rl(f,h) \cong g * h$.
\end{proposition}
\pf
The counit 
$\epsilon: f * g \rightarrow 1: B \rightarrow B$ 
of the adjunction $f \dashv g$ exhibits $g$ as the right lifting
$\Rl(f,1)$
of $1:  B \rightarrow B$ through $f$, and this right 
lifting is absolute. Composing then with any $h: C \rightarrow B$
yields the right lifting $\Rl(f,h)$.
\epf

We shall now consider the more specific case when the 
considered bicategory is {\em closed}, i.e such a bicategory 
admits all right Kan extensions and all right liftings
and is locally complete and cocomplete.
The older terminology for this was {\em biclosed}, but the more 
modern terminology recognises that a bicategory may be 
{\em right-closed} or {\em left-closed}, or {\em both} and to 
be both is to be close.\\

In a closed bicategory $\bicatB$, for any objects $A$
and any arrow $h:C \rightarrow B$, one has the 
natural isomorphisms in $f$, $g$:
\begin{fact}\label{compiso}\end{fact}
\begin{center} 
$\bicatB(A,B)(f, ran_g(h)) \cong \bicatB(C,B)(f * g,h) 
\cong \bicatB(C,A)(g, \Rl(f,h))$
\end{center}  
and thus an adjunction
\begin{fact}\label{genadj}\end{fact}
\begin{center}
$\Rl(-,h) \dashv ran_{-}(h) : \bicatB(A,B) \rightarrow {\bicatB(C,A)}^{op}$ 
\end{center}
Note also that according to \ref{sladj1} and \ref{tradj}
\begin{remark}\label{restoadj}
the composite isomorphism \ref{compiso} in the case $C = B$
and $h = 1_B: B \rightarrow B$,
maps bijectively left adjoints $A \rightarrow B$ 
to right adjoints $B \rightarrow A$.
\end{remark}

We can now turn to the specific case of
the closed bicategory $\VMod$ of modules with domains
and codomains {\em small} categories. For categories
$\A$ and $\B$ not necessarily small, a module 
$\xymatrix{\A \ar[r]|{\circ} & \B}$ 
is a functor $\B^{op} \otimes \A \rightarrow \V$ sending 
$(b,a)$ to $f(b,a)$ and a morphism 
$\alpha: f \rightarrow g: \xymatrix{\A \ar[r]|{\circ} & \B}$
of module is a natural transformation with components
$\alpha_{b,a}: f(b,a) \rightarrow g(b,a)$. So the
modules from $\A$ to $\B$ form an ordinary
category $\VMOD(\A,\B)$ that underlines a
$\V'$-category $\VMOD[\A,\B]$ (for a suitable $\V'$
extending $\V$), which is 
a $\V$-category when $\A$ and $\B$ are small.
Of course a module $\xymatrix{ \A \ar[r]|{\circ} & \B}$ 
can equally be seen as a functor $\A \rightarrow [\B^{op}, \V]$, and 
so (to within equivalence) as a cocontinuous 
functor $\p\A \rightarrow [\B^{op},\V]$.
Therefore when $\A$ and $\B$ are small, 
$\VMOD(\A,\B)$ is equivalent to 
$\Co([\A^{op},\V],[\B^{op},\V]$)  
and the cocontinuous functors here are 
equally the left adjoint ones - see \cite{Kel82}, Theorem 
4.82).
For modules $f: \xymatrix{\A \ar[r]|{\circ} & \B}$
and $g: \xymatrix{\B \ar[r]|{\circ} & \C}$ we have 
a composite $gf: \xymatrix{\A \ar[r]|{\circ} & \C}$
defined by 
\begin{fact}$\;$\label{M1.1}\end{fact}
\begin{center}
$(gf)(c,a) = \int^b  g(c,b) \otimes f(b,a)$
\end{center}
whenever the integral exists - which surely does
when $\B$ is small. This composition is associative
up to coherent isomorphisms, with as 
identities-to-within-isomorphisms the 
$1_{\A}: \xymatrix{\A \ar[r]|{\circ} & \A}$ given
by the hom functor 
$\A(-,-): \A^{op} \otimes \A \rightarrow \V$;   
so that the modules between {\em small} categories
are the arrows of a bicategory $\VMod$.
This $\VMod$ is {\em closed} as it admits 
all right liftings and all right Kan extensions
as follows. 
Given modules $f: \xymatrix{\A \ar[r]|{\circ} & \B}$,
$g: \xymatrix{\B \ar[r]|{\circ} & \C}$
and $h: \xymatrix{\A \ar[r]|{\circ} & \C}$, we consider
new modules $\Sl f,h \Sr: \xymatrix{\B \ar[r]|{\circ} & \C}$ and
$\Cl g,h \Cr: \xymatrix{\A \ar[r]|{\circ} & \B}$ by
\begin{fact}$\;$\label{M1.2}\end{fact}
\begin{center}
$\Sl f,h \Sr(c,b) = \int_a [f(b,a),h(c,a)] $
\end{center}
and 
\begin{fact}$\;$\label{M1.3}\end{fact}
\begin{center}
$\Cl g,h \Cr (b,a) = \int_c [g(c,b),h(c,a)]$
\end{center}
whenever these integrals exist - so that 
$\Sl f,h \Sr$ certainly exists when $\A$ is small,
while $\Cl g, h \Cr$ certainly exists when 
$\C$ is small. When $gf$ and $\Sl f,h \Sr$ exist,
to give a morphism $gf \rightarrow h$ is to
give a family 
$g(c,b) \otimes f(b,a) \rightarrow h(c,a)$
natural (meaning $\V$-natural) in each variable;
which is equally to give a family 
$g(c,b) \rightarrow [ f(b,a),h(c,a) ]$
natural in each variable, and so as morphism   
$g \rightarrow \Sl f, h \Sr$ of modules.
Accordingly so far as the modules in question exist
we have natural bijections
\begin{fact}$\;$\label{M1.4}\end{fact}
\begin{center}
$\VMOD(\A,\C)(gf,h) \cong \VMOD(\B,\C)(g,\Sl f,h \Sr)
\cong \VMOD(\A,\B)(f,\Cl g,h \Cr).$
\end{center}
exhibiting $\Sl f,h \Sr$ as the {\em right Kan extension
of $h$ along $f$} and $\Cl g,h \Cr$ as the {\em right lifting
of $h$ through $g$}.
 
Each functor $T: \A \rightarrow \B$ gives 
rise to modules $T_*: \xymatrix{\A \ar[r]|\circ & \B}$
and $T^*: \xymatrix{\B \ar[r]|{\circ} & A}$, where
\begin{fact}$\;$\label{M1.8}\end{fact}
\begin{center}
$T_*(b,a) = \B(b,Ta)$ and  $T^*(a,b) = \B(Ta,b)$;
\end{center}
and a natural $\rho: T \rightarrow S: \A \rightarrow \B$
gives rise to evident morphisms of modules
\begin{center}
$\rho_*:T_* \rightarrow S_*$ and  $\rho^*: S^* \rightarrow T^*$. 
\end{center}
Since $\rho \mapsto \rho_*$ and 
$\rho \mapsto \rho^*$ are fully faithful by Yoneda, any 
of $\rho$, $\rho_*$ and $\rho^*$ determines the others.
For a functor $T: \A \rightarrow \B$, and modules
$g: \xymatrix{\B \ar[r]|{\circ} & \C}$
and $h: \xymatrix{\D \ar[r]|{\circ} & \B}$
where $\A$, $\B$ and $\C$ are not necessarily small, 
substituting the values \ref{M1.8} into
\ref{M1.1}, \ref{M1.2} and \ref{M1.3} and using the Yoneda 
isomorphisms gives natural isomorphisms
\begin{fact}$\;$\label{M1.10}\end{fact}
\begin{center}
$(gT_*)(c,a) \cong g(c,Ta) \cong \Sl T^*,g \Sr(c,a)$
\end{center}
\begin{fact}\label{M1.11}\end{fact}
\begin{center}
$(T^*h)(a,d) \cong h(Ta,d) \cong \Cl T_*,h \Cr (a,d)$.
\end{center}
In particular for $T: \A \rightarrow \B$ and
$S: \B \rightarrow \C$, one has natural isomorphisms
\begin{center}
 ${(ST)}_* \cong S_* T_*$ and $(ST)^* \cong T^* S^*$;
\end{center} 
so that there are locally fully faithful embeddings of bicategories
\begin{fact}$\;$\label{M1.13}\end{fact}
\begin{center}
${(-)}_*: \VCat \rightarrow \VMod$ $\;\;$
and 
$\;\;$
${(-)}^*: \VCat \rightarrow {\VMod}^{coop}$.
\end{center}

We shall now discuss {\em adjunctions} in $\VMod$.
Consider a functor $T: \A \rightarrow \B$ in $\VCat$ 
and a module $k: \xymatrix{\D \ar[r]|{\circ} & \B}$
in $\VMod$.
Taking the isomorphism \ref{M1.11}, 
$\Cl  T_*, k \Cr \cong T^* k$ shows that
the right lifting of $k$ through $T_*$ is absolute
so by \ref{sladj1}, $T_*$ has right adjoint $T^*$.
In any bicategory left adjoints are often called 
{\em maps} and it is common to denote by $f^*$
the right adjoint of $f$. This is then consistent with
the above notation $T_* \dashv T^*$ above provided
that we identify functor with their associated 
module via the fully faithful $(-)_*$ of
\ref{M1.13}. However the functors are not 
the only maps in $\VMod$ that we shall characterise
now.\\

Consider the following situation in $\VMod$
\begin{center}
$\xymatrix{ & \D \ar[d]|{\circ}^{h} & \\
&  \C \ar[ld]|{\circ}_{g} \ar[rd]|{\circ}^{ \Rl(f, g) = \Cl f, g \Cr } &  \\
\B & &  \A \ar[ll]|{\circ}^{f} \ar@{}[llu]|{\Leftarrow}
}$
\end{center}
where the bottom triangle is a right lifting.
The module $h$ gives by restriction for each object 
$d$ of $\D$ 
a functor $H_d:{\C}^{op} \rightarrow \V$.
Similarly $f$ gives for each $a$ in $\A$
a functor $F_a: {\B}^{op} \rightarrow \V$ and 
the module $g$ corresponds to a functor
$G: {\B}^{op} \rightarrow [\C,\V]$.
One may check straightforwardly that
\begin{fact}\label{rlpre} 
the composition with $h$ respects the 
right lifting $\Cl f,g \Cr$ of $g$ through $f$.
\end{fact} 
if and only if
\begin{fact}\label{colprelim}
for all $a$ in $\A$ and all $d$ in $\D$,
$H_d * - : [\C,\V] \rightarrow \V$
preserves the limit $\{ F_a,G \}$.
\end{fact} 
Writing now $G'$ for the functor $\C \rightarrow [{\B}^{op},\V]$ 
corresponding to $g$, one gets according 
to \ref{fac5.3} 
that the statement above is also equivalent to the fact that
\begin{fact}\label{limprecol} 
for all $a$ in $\A$ and all $d$ in $\D$,
the colimit $H_d * G'$ is preserved
by $\{ F_a, - \}: [{\B}^{op},\V] \rightarrow \V$.
\end{fact} 

At the light of these observations,
we may now reformulate Proposition
\ref{sladj1} for the case when 
the ambient bicategory is $\VMod$.
Note that the assertion $(ii)$ below is the direct translation by 
\ref{rlpre}-\ref{limprecol} of the fact that
for any module $g: \xymatrix{\C \ar[r]|{\circ} & \B}$ the right 
lifting $\Cl f,g \Cr$
is respected by composition with any module 
$\xymatrix{ \I \ar[r]|{\circ} & \C}$.
\begin{proposition}\label{fladj1}
Given small categories $\A$ and $\B$ and
a functor $F: \A \rightarrow [{\B}^{op},\V]$ that corresponds
to the module $f: \xymatrix{\A \ar[r]|{\circ} & \B}$, the following
conditions on $F$ or $f$ are equivalent:
\begin{itemize}
\item $(i)$ $f$ (as a module) has a right adjoint $f^*$;
\item $(ii)$ for each $a$ in $\A$, the representable functor 
$\{Fa,-\}$ $=$ $[{\B}^{op}, \V](Fa,-): [{\B}^{op},\V] \rightarrow \V$ is 
cocontinuous;

\item $(iii)$
for each pair $a$,$c$ in $\A$, the
functor $Fc*-:[\B, \V] \rightarrow \V$
preserves the limit $\{ Fa, Y' \}$ 
of $Y':\B^{op} \rightarrow [\B,\V]$
weighted by $Fa:{\B}^{op} \rightarrow \V$.

\item $(iv)$ for each pair $a$,$c$ in $\A$, the representable
functor $[{\B}^{op},\V](Fa,-):[{\B}^{op}, \V] \rightarrow \V$
preserves the colimit $Fc * Y$ of $Y:\B \rightarrow [{\B}^{op},\V]$
weighted by $Fc:{\B}^{op} \rightarrow \V$ giving (since 
$Fc *Y \cong Fc$ by Yoneda) an isomorphism
$Fc * [{\B}^{op},\V](Fa,Y-) \cong [{\B}^{op},\V](Fa,Fc)$. 
\end{itemize}  
When these are satisfied, we have
\begin{fact}$\;$\label{B2.17}\end{fact}
\begin{center}
$f^*(a,b) \cong [\B^{op},\V](Fa,Yb)$.
\end{center}
\end{proposition}

Proposition \ref{fladj1} gives 
in the particular case $\A = \I$: 
\begin{proposition}\label{fladj2}
Given a small category $\B$ and
a presheaf  $F: {\B}^{op} \rightarrow \V$ that corresponds
to the module $f: \xymatrix{\I \ar[r]|{\circ} & B}$, the following
conditions on $F$ or $f$ are equivalent:
\begin{itemize}
\item $(i)$ $f$ (as a module) has a right adjoint $f^*$;
\item $(ii)$ the representable functor 
$\{F,-\} = [{\B}^{op}, \V](F,-): [{\B}^{op},\V] \rightarrow \V$ is 
cocontinuous;
\item $(iii)$ the functor $F*-:[\B, \V] \rightarrow \V$
preserves the limit $\{ F, Y' \}$ 
of $Y':\B^{op} \rightarrow [\B,\V]$
weighted by $F:{\B}^{op} \rightarrow \V$.

\item $(iv)$ the representable
functor $\{F,-\} = [{\B}^{op},\V](F,-):[{\B}^{op}, \V] \rightarrow \V$
preserves the colimit $F * Y$ of $Y:\B \rightarrow [{\B}^{op},\V]$
weighted by $F:{\B}^{op} \rightarrow \V$ giving (since 
$F *Y \cong F$ by Yoneda) an isomorphism
$F * \{F,Y-\} \cong \{F,F\}$; 
\end{itemize}  
When these are satisfied, we have
\begin{fact}$\;$\label{B2.17bis}\end{fact}
\begin{center}
$f^*(b) \cong [\B^{op},\V](F,Yb)$.
\end{center}
\end{proposition}

Once again there are dual statements to \ref{fladj1} and 
\ref{fladj2} obtained by translating part of \ref{tradj}. 
In detail:
\begin{proposition}
For any functor $G: {\A}^{op} \rightarrow [\B,\V]$ that corresponds
to the module $g: \xymatrix{\B \ar[r]|{\circ} & \A}$, the following
conditions on $G$ or $g$ are equivalent:
\begin{itemize}
\item $(i)$ $g$ (as a module) has a left adjoint $g_*$;
\item $(ii)$ for each $a$, the representable functor 
$\{Ga, -\} = [\B, \V](Ga,-): [\B,\V] \rightarrow \V$ is 
cocontinuous;
\item $(iii)$ for each pair $a$,$c$ in $\A$, the 
functor $Gc * -: [\B^{op}, \V] \rightarrow \V$
preserves the limit $\{Ga, Y \}$ of 
$Y: \B \rightarrow [{\B}^{op},\V]$ weighted by 
$Ga: \B \rightarrow \V$;
\item $(iv)$ for each pair $a$,$c$ in $\A$, the representable
functor $[\B,\V](Ga,-):[\B, \V] \rightarrow \V$
preserves the colimit $Gc * Y'$ where $Y':{\B}^{op} \rightarrow
[\B,\V]$.
\end{itemize}  
When these are satisfied, we have
\begin{fact}$\;$\label{B2.18}\end{fact}
\begin{center}
$g_*(a,b) \cong [\B,\V](Ga,Y'b)$.
\end{center}
\end{proposition}

Eventually for the case $\A = \I$ one has
\begin{proposition}\label{gradj}
For any weight $G: \B \rightarrow \V$ that corresponds
to the module $g: \xymatrix{\B \ar[r]|{\circ} & \I}$, the following
conditions on $G$ or $g$ are equivalent:
\begin{itemize}
\item $(i)$ $g$ (as a module) has a left adjoint $g_*$;
\item $(ii)$ the representable functor 
$\{G,-\}= [\B, \V](G,-): [\B,\V] \rightarrow \V$ is 
cocontinuous;
\item $(iii)$ the 
functor $G * -: [\B^{op}, \V] \rightarrow \V$
preserves the limit $\{G, Y \}$ of 
$Y: \B \rightarrow [{\B}^{op},\V]$ weighted by 
$G: \B \rightarrow \V$. 
\item $(iv)$ the representable
functor $\{G,-\}:[\B, \V] \rightarrow \V$
preserves the colimit $G * Y'$.
\end{itemize}  
When these are satisfied, we have
\begin{fact}$\;$\label{B2.18bis}\end{fact}
\begin{center}
$g_*(b) \cong [\B,\V](G,Y'b)$.
\end{center}
\end{proposition}

For any presheaves $F,H: \B^{op} \rightarrow \V$
on a small $\B$, that
correspond respectively to modules 
$f,h: \xymatrix{\I \ar[r]|{\circ} & \B}$,
$[\B^{op}, \V](F,H)$ is the right lifting of $h$ through 
$f$. Therefore as a consequence of Proposition \ref{lemadj}: 
\begin{proposition}
\label{Homladj}
Given a module $f: \xymatrix{\I \ar[r]|{\circ} & \B}$
corresponding to a presheaf $F: \B^{op} \rightarrow \V$,
if $f$ has a right adjoint $g: \xymatrix{\B \ar[r]|{\circ} & \I}$
then for any left module $h: \xymatrix{\I \ar[r]|{\circ} & \B}$ 
corresponding to a presheaf $H: \B^{op} \rightarrow \V$
\begin{center}
$[\B^{op}, \V](F,H) \cong g * h$.
\end{center}
\end{proposition}

\begin{proposition}\label{yac}
For any $G: \B \rightarrow \V$ with $\B$ small,
that corresponds to a  
module $g: \xymatrix{ \B \ar[r]|{\circ} & \I }$, 
the following assertions are equivalent:
\begin{itemize} 
\item $(i)$ $g$ is a right adjoint;\\
\item $(ii)$ $-* G: [{\B}^{op}, \V] \rightarrow \V$
is representable;\\
\item $(iii)$ $-* G: [{\B}^{op}, \V] \rightarrow \V$
is continuous.
\end{itemize}
\end{proposition}
\pf
$(i) \Rightarrow (ii)$. If $g$ has left adjoint 
$f: \xymatrix{ I \ar[r]|{\circ} & \V}$ that corresponds to
a functor $F: \B^{op} \rightarrow \V$ then
$\{F,-\} \cong -*G \cong G*-$ by \ref{Homladj}. 
$G*-$ being representable is therefore continuous.\\
$(ii) \Rightarrow (iii)$ well-known.\\
$(iii) \Rightarrow (i)$.
$G*- = -*G: [{\B}^{op},\V] \rightarrow \V$
being continuous preserves in particular the limit
$\{ G,Y\}$ of $Y: \B \rightarrow [\B^{op},\V]$
weighted by $G$ and by Proposition \ref{gradj} 
this equally asserts that $g$ is a right adjoint.
\epf

In \cite{Str83} Street characterised 
right adjoints modules $g: \xymatrix{ \B \ar[r]|{\circ} &  \I}$
with $\B$ small
as the ones for which {\em the corresponding
presheaf $G: \B \rightarrow \V$ is a weight
of absolute colimits}. 
A similar argument to the one used for the proof 
$(i) \Rightarrow (ii)$ above, yields half of the result. 
Actually if $G$ is a weight of absolute colimits 
then the colimit $G*Y'$ is certainly preserved 
by $\{G,-\}$ and then by \ref{gradj}, $g$ is right adjoint.\\

An object $a$ of a category $\A$ (especially of a cocomplete
$\A$) is said in \cite{Kel82} to be {\em small projective}
when the representable $\A(a,-)$ preserves all small colimits.
Of course the representables in $[\B^{op},\V]$ are among 
the small projectives: for  $[\B^{op},\V](Yb,-)$ is isomorphic
to the evaluation $E_b$ which preserves all colimits.   
Proposition \ref{fladj1} above asserts in particular 
that $f: \xymatrix{\A \ar[r]|{\circ} & \B}$ has a right 
adjoint (or is a map) precisely when each 
$f(-,a): \B^{op} \rightarrow \V$ is a small projective 
in $[\B^{op},\V]$. Note the special case of a functor 
$T: \A \rightarrow \B$ where $T_*(-,a) = \B(-,Ta)$
and the representable $\B(-,Ta)$ is a small projective.
There is another way of looking at the small projectives
in the category $[\B^{op},\V]$ when $\B$ is small:
an object $\phi$ of $[\B^{op},\V]$ is a weight,
and $[\B^{op},\V](\phi,\psi)$ is just the 
$\phi$-weighted limit of $\psi: \B^{op}\rightarrow \V$.
To say, therefore that $\phi$ is a small projective
is equally to say that $\phi$-limits commute
in $\V$ with all small colimits. As mentioned earlier
we shall write $\Q$ for the class of small projective
weights. In the language of the previous section, 
$\Q = \p^-$, and according to \ref{yac} also $\Q = \p^+$.
For any small $\B$, as $\Q$ is saturated, the closure 
$\Q(\B)$ of $\B$ in $[\B^{op},\V]$ under $\Q$-colimits
is $\Q[\B]$, as $\Q$ is equally the class 
of weights corresponding
in the sense of \ref{fladj1} to left modules,  
$\Q(\B)$ is the {\em Cauchy-completion} of $\B$
in the sense of \cite{Law73}.\\

When the category $\B$ is not small, we can still
consider the closure $\Q(\B)$ of $\B$ in $[\B^{op},\V]$
(or equally in $\p(\B)$) under $\Q$-colimits. On the other
hand, we can consider the full subcategory of $\p(\B)$
given by its small projectives that is $\p(\B)_\p$.
By Proposition \ref{pro5.2}, $\p(\B)_{\p}$ is closed
in $\p(\B)$ under $\Q$-colimits. Since the representables
lie in $\p(\B)_{\p}$, it follows that 
$\Q(\B) \subseteq {\p(\B)}_{\p}$. When $\V = Set$ the converse
is true, so that $\Q(\B)$ consists even for large $\B$
of the small projectives in $\p(\B)$. To check this 
last point one can adapt the proof of  \cite{Kel82} to show
that any accessible $F: \B^{op} \rightarrow Set$ such that 
$[\B^{op},Set](F,-)$ preserves small colimits is a retract
of representables, thus a $\Q=\p^{+}$-weighted colimit of 
representables, and thus in $\Q(\B)$.\\

For $\B$ small, the operations \ref{B2.17bis} and 
\ref{B2.18bis} above define actually two functors: 
\begin{fact}$\;$\label{B3.1}\end{fact}
\begin{center}
$[\B^{op},\V](-,Y): [{\B}^{op},\V] \rightarrow  {[\B,\V]}^{op}$ 
\end{center}
and 
\begin{fact}$\;$\label{B3.2}\end{fact}
\begin{center}
$[\B,\V](-,Y'): [\B,\V] \rightarrow  {[\B^{op},\V]}^{op}$.
\end{center}
Letting $\bicatB = \VMod$,
$A =  \I$ and $h:C \rightarrow B$ 
$=$ $1_{\B}: \B \rightarrow \B$ in 
\ref{genadj} we get that their underlying Set-functors
define an adjoint pair:
\begin{fact}$\;$\label{Isbadj}\end{fact}
\begin{center} 
$\VCat(\B^{op},\V) \rightharpoonup {\VCat({\B},\V)}^{op}$.
\end{center}
The adjunction above goes back in the case $\V=Set$ to
Isbell. Let us call it the {\em Isbell adjunction}.
Also by Remark \ref{restoadj}, Propositions \ref{fladj2} and \ref{gradj},
\begin{proposition} 
The Isbell adjunction as in \ref{Isbadj} restricts to an 
equivalence  
\begin{center}
${\Q( \B^{op} ) }^{op} \cong \Q(\B)$
\end{center}
where $\Q(\B)$ consists of the small projectives in
$[{\B}^{op},\V]$.
\end{proposition}
We shall further study more in detail the Isbell adjunction
and shows that it lifts to a $\V$-adjunction.\\

Recall that for a small $\B$, the functor category $[\B^{op},\V]$
is the free cocomplete category on $\B$, in the sense that,
for each cocomplete $\C$, we have an adjoint equivalence
of categories:
\begin{fact}$\;$\label{B3.3}\end{fact}
\begin{center}
$\xymatrix{ \Co[[\B^{op},\V], \C] \ar@<1ex>[r]^(.65){[Y,1]} & 
[\B,\C] \ar@<1ex>[l]^(.35){Lan_Y}}$
\end{center} 
where $[Y,1]$ sends the cocontinuous 
$S:[{\B}^{op},\V] \rightarrow \C$ to its restriction
$SY: \B \rightarrow \C$, while $Lan_Y$ sends $G:\B \rightarrow \C$
to $Lan_Y(G) = -*G:[\B^{op},\V] \rightarrow \C$ which in fact 
has the right adjoint 
$\tilde{G}: \C \rightarrow [\B^{op},\V]$ where
\begin{fact}$\;$\label{B3.4}\end{fact}
\begin{center}
$\tilde{G}c = \C(G-,c)$
\end{center}
as in the diagram
\begin{fact}$\;$\label{B3.5}\end{fact}
\begin{center}
$\xymatrix{\C \ar@<1ex>[rr]^{\tilde{G}} 
& 
& [\B^{op},\V] \ar@<1ex>[ll]^{-*G}\\
\\
& \B \ar[luu]^G \ar[ruu]_Y & 
}$
\end{center}
Recall too that there is a dual version involving
$\hat{G}: \C^{op} \rightarrow [\B,\V]$ where
\begin{fact}$\;$\label{B3.6}\end{fact}
\begin{center}
$\hat{G}c = \C(c,G-)$.
\end{center} 

Now consider the case where $G$ is the Yoneda embedding
\begin{fact}
$y = {Y'}^{op}: \B \rightarrow [ \B, \V]^{op}$
\end{fact}
whereupon \ref{B3.5} becomes
\begin{fact}$\;$\label{B3.8}\end{fact}
\begin{center}
$\xymatrix{[\B,\V]^{op} \ar@<1ex>[rr]^{\tilde{y}} 
& 
& [\B^{op},\V] \ar@<1ex>[ll]^{-*y}\\
\\
& \B \ar[luu]^y \ar[ruu]_Y & 
}$
\end{center}
Here 
\begin{fact}$\;$\label{B3.9}\end{fact}
\begin{center}
$\tilde{y}(\varphi) = [\B,\V]^{op}(y-,\varphi)
= [\B,\V](\varphi,Y'-)$, 
\end{center}
so that 
\begin{fact}$\;$\label{B3.10}\end{fact}
\begin{center}
$\tilde{y}(\varphi) = \psi$ where $\psi(b) = [\B,\V](\varphi, Y'b)$,
\end{center}
as in \ref{B3.2}. Again, since $\psi * y$ is the colimit
in $[\B,\V]^{op}$ of $y: \B \rightarrow {[\B,\V]}^{op}$
weighted by $\psi: \B^{op} \rightarrow \V$, it is equally
the limit $\{\psi, Y'\}$ in $[\B, \V]$ of
$y^{op} = Y': \B^{op} \rightarrow [\B,\V]$ weighted by
$\psi: \B^{op} \rightarrow \V$; and since limits in 
$[\B,\V]$ are formed pointwise from the limits in $\V$,
its value $(\psi * y)b = \{ \psi, Y' \}b$ at $b \in \B$
is given by
\begin{fact}$\;$\label{B3.11}\end{fact}
\begin{center}
$(\psi * y)(b) = \{\psi, Yb \} = [{\B}^{op},\V](\psi, Yb)$.
\end{center}  
An immediate consequence of this, using \ref{B3.6}, is the 
observation that 
\begin{fact}$\;$\label{B3.12}\end{fact}
\begin{center}
$-*y = \hat{Y}^{op}: [\B^{op},\V] \rightarrow [\B,\V]^{op}$
\end{center}
which may be combined with the trivially true
\begin{fact}$\;$\label{B3.13}\end{fact}
\begin{center}
$\tilde{y} 
= {({Y'}^{op})}^{\sim}:[ \B,\V]^{op} \rightarrow [\B^{op},\V]$
\end{center}
in order to write the adjunction in \ref{B3.8} as 
\begin{fact}$\;$\label{B3.14}\end{fact}
\begin{center}
${\hat{Y}}^{op} \dashv {({Y'}^{op})}^{\sim}$
\end{center}
Omitting the $b$ in \ref{B3.11} gives us a form which we may compose
with \ref{B3.9}, namely
\begin{fact}$\;$\label{B3.15}\end{fact}
\begin{center}
$\psi * y = [\B^{op},\V](\psi,Y-)$
\end{center}
so that - compare \ref{B3.10} - we have
\begin{fact}$\;$\label{B3.16}\end{fact}
\begin{center}
$\psi * y = \varphi$ where $\varphi(b) = [\B^{op}, \V](\psi, Yb)$,
\end{center} 
as in \ref{B3.1}.
Now we observe that there is an alternative 
analysis of this equivalence between the 
$\psi$'s in $\Q(\B)$ and the $\varphi$'s in $\Q(\B^{op})$.
Accordingly to the general picture in \ref{B3.5}
referring to $[\B^{op},\V]$ as the free cocompletion
of $\B$, when $\psi \in \Q(\B)$, the cocontinuous
$[\B^{op},\V](\psi,-):[\B^{op},V] \rightarrow \V$
is given by:
\begin{fact}\label{B4.5}\end{fact}
\begin{center}
$\xymatrix{\V \ar@<1ex>[rr]^{\tilde{\varphi}} 
& 
& [\B^{op},\V] \ar@<1ex>[ll]^{\{ \psi,- \} =-*\varphi}\\
\\
& \B \ar[luu]^{\varphi = [\B^{op},\V](\psi,Y-)} \ar[ruu]_Y & 
}$
\end{center}
In more detail, $\psi$ lies in $\Q(\B)$ precisely when
$[\B^{op},\V](\psi,-)$ is cocontinuous (or equally
left adjoint); and then $[\B^{op},\V](\psi,-) = -* \varphi$
where $\varphi = [\B^{op},\V](\psi,Y-) = \psi*y$.
On the other hand, if we start with $\varphi:\B \rightarrow \V$
and ask whether it is of the form $[\B^{op},\V](\psi,Y-)$
for some $\psi \in \Q(\B)$, we observe that it so 
precisely when 
$-* \varphi: [\B^{op},\V] \rightarrow \V$
is {\em representable}.

\end{section}

\bibliographystyle{alpha}

\begin{thebibliography}{}

\bibitem[ABLR02]{ABLR02}
{\sc Jiri Adamek, Francis Borceux, Stephen Lack and Jiri Rosicky},\\
{\it A classification of accessible categories},\\
Journal of Pure and Applied Algebra 175, issues 1-3, 2002, 7-30
\bibitem[AK88]{AK88}
{\sc M.H. Albert, G.M. Kelly},\\
{\it The closure of a class of colimits},\\
Journal of Pure and Applied Algebra 51, 1988, 1-17
\bibitem[AR94]{AR94}
{\sc Jiri Adamek and Jiri Rosicky},\\
{\it Locally presentable and accessible categories},\\
Cambridge University Press, 94.
\bibitem[Bet85]{Bet85}
{\sc R. Betti},\\
{\it Cocompleteness over coverings},\\
J. Austral. Math. Soc. (Series A) 39, 1985, 169-177.
\bibitem[BCSW83]{BCSW83}
{\sc R. Betti, A. Carboni, R. Street, R. Walters},\\
{\it Variation through enrichment },\\
J. Pure and App. Algebra 29, 1983, 109-127.
\bibitem[BQR98]{BQR98}
{\sc F. Borceux, C. Quintero, J. Rosicky},\\
{\it A theory of enriched sketches},\\
Theory and Applications of Categories, Vol. 4, No. 3, 1998, 47-72.
\bibitem[DaSt97]{DaSt97}
{\sc B. Day, R. Street}
{\it Monoidal bicategories and Hopf algebroids},\\
Advance in maths. 129, 1997, 99-157.
\bibitem[KeLa80]{KeLa80}
{\sc G.M. Kelly, Laplaza},\\
{\it Compact closed categories},\\
J. Pure and App. Algebra 19, 1980, 193-213.    
\bibitem[Kel82]{Kel82}
{\sc G.M. Kelly},\\
{\it Basic concepts of enriched category theory},\\
London Mathematical Society Lecture Note Series 64
- Cambridge University Press 82.
\bibitem[Kel82-2]{Kel82-2}
{\sc G.M. Kelly},\\
{\it Structures defined by finite limits in the enriched
context, 1},\\
Cahiers de Top. et Geo. diff. Vol. XXIII-1, 1982.
\bibitem[Law73]{Law73}
{\sc F.W. Lawvere},\\
{\it Metric spaces, generalized logic, and closed categories},\\
Rend. Sem. Mat. Fis. di Milano 43, 73, 135-166,\\
Reprints in Theory and Applications of Categories, No. 1, 2002, 1-37.
\bibitem[MP89]{MP89}
{\sc M. Makkai and R. Par\'e},\\
{\it Accessible categories: the foundations of categorical 
model theory},\\
Contemporary Math. 104, AMS, Providence, 1989.
\bibitem[Str83]{Str83}
{\sc R. Street},\\
{\it Absolute colimits in enriched categories},\\
Cahiers de Top. et Geo. diff. Vol. XXIV-4, 1983.
\bibitem[StWa78]{StWa78}
{\sc R. Street, R. Walters},\\
{\it Yoneda structures on 2-categories},\\
Journal of Algebra 50, 350-379, 1978.
\end{thebibliography}

\end{document}